\documentclass[letterpaper, 10 pt, conference]{ieeeconf}  
         
\usepackage{latexsym,amsmath}

\usepackage{amsfonts}

\def\half{\mbox{$\frac{1}{2}$}}

\newcommand{\shoot}{\mathcal{S}}

\newcommand{\lam}{[\lambda]}
\newcommand{\lamh}{[\hat\lambda]}

\newtheorem{theorem}{Theorem}[section]
\newtheorem{lemma}{Lemma}[section]
\newtheorem{proposition}[theorem]{Proposition}

\newtheorem{assumption}{Assumption}
\newtheorem{remark}{Remark}[section]
\newtheorem{definition}{Definition}[section]

\def\ds{\displaystyle}

\def\bega{\begin{array}}
\def\enda{\end{array}}

\def\bepmatrix{\begin{pmatrix}}
\def\enpmatrix{\end{pmatrix}}

\def\bel{\begin{equation}\label}
\def\eeq{\end{equation}}
\newcommand\ee{\end{equation}}
\def\benl{\begin{equation*}}
\def\eenl{\end{equation*}}

\def\be{\begin{equation}}
\def\beq{\begin{equation}}
\def\bel{\begin{equation}\label}
\def\eeq{\end{equation}}

\newcommand\ba{\begin{array}}
\newcommand\ea{\end{array}}

\def\begi{\begin{itemize}}
\def\endi{\end{itemize}}

\newcommand{\cR}{\mathbb{R}}

\newcommand{\tras}{^\top}
\newcommand{\intT}{\int_0^T }
\newcommand{\ddt}{\frac{\rm d}{{\rm d} t} }
\newcommand{\dtt}{\mathrm{d}t}
\newcommand{\mr}{\mathrm}

\def\C{\mathcal{C}}

\def\J{\mathcal{J}}

\def\L{\mathcal{L}}

\def\P{\mathcal{P}}

\def\U{\mathcal{U}}
\def\V{\mathcal{V}}
\def\W{\mathcal{W}}
\def\X{\mathcal{X}}

\def\hb{\bar{h}}

\def\ub{\bar{u}}
\def\vb{\bar{v}}
\def\wb{\bar{w}}
\def\xb{\bar{x}}
\def\yb{\bar{y}}

\def\xib{\bar{\xi}}

\def\fh{\hat{f}}

\def\ph{\hat{p}}

\def\uh{\hat{u}}
\def\vh{\hat{v}}
\def\wh{\hat{w}}
\def\xh{\hat{x}}

\def\Ih{\hat{I}}
\def\Th{\hat{T}}

\def\Wh{\hat{W}}

\def\eps{\varepsilon}

\IEEEoverridecommandlockouts                              
\title{\LARGE \bf
Convergence of the shooting algorithm \\ for singular  optimal control problems*
}

\author{M. Soledad Aronna$^{1}$ 
\thanks{*This article will appear in the Proceedings of the European Control Conference to be held in Zurich, Switzerland, 2013. 
}
\thanks{$^{1}$M.S. Aronna is a Marie Curie Fellow within the Network SADCO at the Department of Electrical and Electronic Engineering, Imperial College of London,
South Kensington campus,
London, SW7 2AZ,
UK,
        {\tt\small aronna@cmap.polytechnique.fr}}
}

\begin{document}

\maketitle
\thispagestyle{empty}
\pagestyle{empty}

\begin{abstract}

In this article we propose a shooting algorithm for optimal control problems governed by systems that are affine in one part of the control variable. Finitely many equality  constraints on the initial and final state are considered. We recall a second order sufficient condition for weak optimality, and show that it guarantees
 the local quadratic convergence of the algorithm.
We show an example and solve it numerically.
\end{abstract}


\section{INTRODUCTION}\label{SecInt}

We investigate optimal control problems governed by ordinary differential equations that are affine in one part of the control variable. This class of system includes both the totally affine and the nonlinear cases. 
This study is motivated by many models that are found in practice. Among them we can cite the followings: the Goddard's problem analyzed in Martinon et al. \cite{BLMT09,BryHo,Mau76}, other models concerning the motion of a rocket in Lawden \cite{Law63}, Bell and Jacobson \cite{BelJac}, Goh \cite{Goh08}, Oberle \cite{Obe77}, and an optimal production process in Cho et al. \cite{ChoAbaPar93}.

We can find shooting-like methods applied to the numerical solution of partially affine problems in, for instance,  Oberle \cite{OberleThesis,Obe90} and Oberle-Taubert \cite{ObeTau97}, where the authors use a generalization of the algorithm that Maurer \cite{Mau76} suggested for totally affine systems.
These works present interesting implementations of a shooting-like algorithm, but they do not link the convergence of the method with sufficient conditions of optimality as it is done in this article.

In this paper we propose a shooting algorithm which can be also used to solve problems with bounds on the controls. We give a theoretical support to this method, by showing that a second order sufficient condition for optimality proved in Aronna \cite{Aro11} ensures the local quadratic convergence of the algorithm. 

The article is organized as follows. In Section \ref{SecPb} we give the statement of the problem,  the main definitions and assumptions, and a first order optimality condition. The shooting algorithm is described in Section \ref{SecShoot}. In Section \ref{SecSSC} we recall a second order sufficient condition for weak optimality. We state the main result of the article in Section \ref{SecConv}. In Section \ref{SecEx} we work out an example and solve it numerically. 

\vspace{6pt}

\noindent{NOTATIONS.}
Let $h_t$ denote the value of function $h$ at time $t$ if $h$ is a function that depends only on $t,$ and $h_{i,t}$ the $i$th component of $h$ evaluated at $t.$ Partial derivatives of a function $h$ of $(t,x)$ are referred as $D_th$ or $\dot{h}$ for the derivative in time, and $D_xh$ or $h_x$ for the differentiations with respect to space variables. The same convention  is extended to higher order derivatives.
By $L^p(0,T;\cR^k)$ we mean the Lebesgue space with domain equal to the interval $[0,T]\subset \cR$ and with values in $\cR^k.$ The notation $W^{q,s}(0,T;\cR^k)$ refers to the Sobolev spaces.

\section{Statement of the Problem}\label{SecPb}

We study the optimal control
problem (P) given by
\begin{align}
&\label{cost} J:=\varphi_0(x_0,x_T)\rightarrow
\min,\\
&\label{stateeq}\dot{x}_t= F(x_t,u_t,v_t)=\sum_{i=0}^m v_{i,t}
f_i(x_t,u_t) ,\  {\rm  a.e.}\ {\rm  on}\ [0,T],\\
& \label{finaleq} \eta_j(x_0,x_T)=0,\quad
\mathrm{for}\
j=1\hdots,d_{\eta}.
\end{align}
Here $f_i:\cR^{n+l}\rightarrow
\cR^n$ for $i=0,\hdots,m,$ 
$\varphi_0:\cR^{2n}\rightarrow \cR,$ 
$\eta_j:\cR^{2n}\rightarrow \cR$ for
$j=1,\hdots,d_{\eta}$ and we put, in sake of simplicity of notation, $v_0\equiv 1$ which is not a variable. 
The \textit{nonlinear control} $u$ belongs to $\U:=L^{\infty}(0,T;\cR^l),$ while 
$\V:=L^{\infty}(0,T;\mathbb{R}^m)$ denotes the space of {\it affine controls} $v,$ and 
$\mathcal{X}:=W^{1,\infty}(0,T;\mathbb{R}^n)$ refers to the state space.
 When needed, we write $w=(x,u,v)$ for a point in
$\mathcal{W}:=\mathcal{X}\times \U \times \V.$ 
Assume throughout the article that data functions $\varphi_0,$ $f_i$ and $\eta_j$ have Lipschitz-continuous second derivatives.
A \textit{trajectory} is an element $w\in\mathcal{W}$
that satisfies the state equation \eqref{stateeq}.
If in addition, the constraints in \eqref{finaleq} 
hold, we say that $w$ is a \textit{feasible trajectory}
of problem (P).

Set $\X_*:= W^{1,\infty}(0,T;\cR^{n,*})$ the space of Lipschitz-continuous  functions with values in the $n-$dimensional space of row-vectors with real components $\cR^{n,*}.$
Consider an element $\lambda:=(\beta,p)\in \cR^{d_{\eta},*}\times \X_*$ and define
the \textit{pre-Hamiltonian} function
\be\label{Hdef}
H\lam(x,u,v,t):=p_tF(x,u,v),
\ee
the \textit{initial-final Lagrangian} function
\be
\label{elldef}
\ell[\lambda](\zeta_0,\zeta_T)
:=\varphi_0(\zeta_0,\zeta_T)+\sum_{j=1}^{d_{\eta}}\beta_j \eta_j(\zeta_0,\zeta_T),
\ee
and the \textit{Lagrangian} function
\benl
\L\lam (w):=\ell\lam(x_0,x_T)
+\intT p_t\left(F(x_t,u_t,v_t)-\dot x_t\right)\dtt.
\eenl

Throughout the article we study a nominal feasible trajectory $\wh=(\xh,\uh,\vh),$ that we assume to be smooth.
We present now  an hypothesis for the endpoint constraints. Consider the mapping
\benl 
\ba{rcl}
G\colon  \cR^n\times \U \times \V &\rightarrow &\cR^{d_{\eta}}\\
 (x_0,u,v) & \mapsto &\eta(x_0,x_T),
\ea
\eenl
where $x_t$ is the solution of \eqref{stateeq} associated with $(x_0,u,v).$

\begin{assumption}
\label{lambdaunique}
The derivative of $G$ at $(\xh_0,\uh,\vh)$ is onto.
\end{assumption}

The latter hypothesis is usually known as \textit{qualification of the endpoint equality constraints.} 

\begin{definition}
It is said that the  feasible trajectory $\wh$ is a \textit{weak minimum} of problem (P) if there exists $\varepsilon>0$ such that $\wh$ is a minimum in the set of feasible trajectories $w=(x,u,v)$ satisfying 
\benl
\|x-\xh\|_{\infty}<\varepsilon,\quad  \|u-\uh\|_{\infty}<\varepsilon\quad  \|v-\vh\|_{\infty}<\varepsilon.
\eenl
\end{definition}

The following first order necessary condition holds for $\wh.$
See the book by Pontryagin et al. \cite{PBGM} for a proof.
\begin{theorem}
\label{Mult}
 Let $\wh$ be a weak solution satisfying Assumption \ref{lambdaunique}, then there exists a unique $\hat\lambda=(\hat\beta,\ph) \in \cR^{d_{\eta},*} \times \X_*$ such that $\ph$ is solution of the \textit{costate equation} 
\be\label{costateeq}
-\dot{\ph}_t =D_xH[\hat\lambda](\xh_t,\uh_t,\vh_t,t),\quad {\rm  a.e.}\ {\rm on}\ [0,T],
\ee
with \textit{transversality conditions} 
\begin{align}
 \label{p0} \ph_0 & =-D_{x_0}\ell[\hat\lambda](\xh_0,\xh_T),\\
 \label{pT} \ph_T & =D_{x_T}\ell[\hat\lambda](\xh_0,\xh_T),
\end{align}
and the \textit{stationarity condition}
\be
\label{stationarity}
\left\{
\ba{l}
\vspace{3pt} \ds H_u\lamh(\xh_t,\uh_t,\vh_t,t)=0,\\
H_v \lamh(\xh_t,\uh_t,\vh_t,t)=0,
\ea
\right.
\quad {\rm a.e.}\ {\rm on}\ [0,T],
\ee
is verified.
\end{theorem}

Throughout this article $\wh$ is considered to be a weak solution and thus, it satisfies \eqref{stationarity} for its unique associated multiplier $\hat\lambda.$ Furthermore, note that since $v$ appears linearly in $H$ we have that $D^2_{(u,v)^2}H\lamh(\xh_t,\uh_t,\vh_t,t)$ is a singular matrix on  $[0,T].$ Therefore, $\wh$ is a {\it singular solution} (as defined in \cite{BelJac} and \cite{BryHo}).

\section{The shooting algorithm}\label{SecShoot}

The purpose of this section is to present an appropriate numerical scheme to solve the problem (P).
More precisely, we investigate the formulation and the convergence of an algorithm that approximates an optimal solution provided an initial estimate exists.

\subsection{Optimality system} In what follows we use the first order optimality conditions \eqref{stationarity} to provide a set of equations from which we can determine $\wh.$ We obtain an optimality system in the form of a {\it two-point boundary value problem} (TPBVP). 

Throughout the rest of the article we assume, in sake of simplicity, that whenever some argument of $f_i,$ $H,$ $\ell,$ $\L$ or their derivatives is omitted,  they are evaluated at $\wh$ and $\hat\lambda.$

We shall recall that for the case where all the control variables appear nonlinearly ($m=0$), the classical technique is using the stationarity equation 
\be 
\label{Hu0}
H_u\lamh (\wh)=0,
\ee
to write $\uh$ as a function of $(\xh,\hat\lambda).$ This procedure is also detailed in \cite{MauGil75} and \cite{Tre12}. One is able to do this by assuming, for instance, the \textit{strengthened Legendre-Clebsch condition} 
\be
\label{strongLC}
H_{uu}\lamh (\wh) \succ 0.
\ee
In case \eqref{strongLC} holds, due to the Implicit Function Theorem, we can write $\uh=U\lamh(\xh)$ with $U$ being a smooth function. Hence, replacing the occurrences of $\uh$  by $U\lamh(\xh)$ in the  state and costate equations yields a two-point boundary value problem.

On the other hand, when the system is affine in all the control variables ($l=0$), we cannot eliminate the control from the equation $H_v=0$ and, therefore, a different technique is employed (see e.g. \cite{Mau76,ABM11,Tre12}).  
The idea is to consider an index $1\leq i \leq m,$ and to take ${{\rm d}^{M_i}H_v}/{{\rm d}t^{M_i}}$ to be the lowest order derivative of $H_v$ in which $\vh_i$ appears with a coefficient that is not identically zero.
Goh \cite{Goh66a,GohThesis}, Kelley et al. \cite{KelKopMoy67} and Robbins \cite{Rob67} proved that $M_i$ is even.
This implies that the control does not appear the first time we derive  $H_v$ with respect to time, i.e. $\dot H_v$ depends only on $\xh$ and $\hat\lambda$ and consequently, it is differentiable in time. Thus the expression
\be
\label{ddotHv0}
\ddot{H}_v \lamh (\wh)=0
\ee
is well-defined. The control $\vh$ can be retrieved from \eqref{ddotHv0}  provided that, for instance, the \textit{strengthened generalized Legendre-Clebsch condition} 
\be
\label{strongenLC}
-\frac{\partial \ddot H_v}{\partial v}\lamh (\wh) \succ 0
\ee
holds (see Goh \cite{GohThesis,Goh95,Goh08}). 
In this case, we can write $\vh =  V\lamh(\xh)$ with $V$ being differentiable.
By replacing $\vh$ by $V\lamh(\xh)$ in the state-costate equations, we get an optimality system in the form of a boundary value problem.

In the problem studied here, where $l> 0$ and $m> 0,$ we aim to use both equations \eqref{Hu0} and \eqref{ddotHv0} to retrieve the control $(\uh,\vh)$ as a function of  the state $\xh$ and the multiplier $\hat\lambda.$
We next describe a procedure to achieve this elimination
that was proposed in Goh \cite{Goh95,Goh08}. 
First let as recall a necessary condition proved in Goh \cite{Goh66} and in \cite[Lemma 3.10 and Corollary 5.2]{Aro11}. Define $[f_i,f_j]^x:=(D_xf_j)f_i-(D_xf_i)f_j,$ which is referred as the {\it Lie bracket} in the variable $x$ of $f_i$ and $f_j.$ 
\begin{lemma}[Necessary conditions for weak optimality]
\label{NC}
If $\wh$ is a smooth weak minimum for (P) satisfying Assumption \ref{lambdaunique}, then
\begin{gather}
\label{Huv} H_{uv}\equiv 0,\\
\label{Lie}
\hat{p}[f_i,f_j]^x=0,\quad \text{for}\ i,j=1,\dots,m.
\end{gather}
\end{lemma}

Let us show that $H_v$ can be differentiated twice with respect to the time variable, as it was done in the totally affine case.
Observe that \eqref{Hu0} may be used to write $\dot\uh$ as a function of $(\hat\lambda,\wh).$ 
In fact, in view of Lemma \ref{NC},
the coefficient of $\dot\vh$ in $\dot{H}_u$ is zero. Consequently,
\be 
\label{dotHu0}
\dot H_u=\dot H_u\lamh (\xh,\uh,\vh,\dot\uh)=0
\ee
and, if the strengthened Legendre-Clebsch condition \eqref{strongLC} holds, $\dot{\uh}$ can be eliminated from \eqref{dotHu0} yielding
\be 
\label{dotu}
\dot\uh = \Gamma\lamh(\xh,\uh,\vh).
\ee
Take now an index $i=1,\dots,m$ and observe that
\benl
0=\dot{H}_{v_i} = \ddt \,\ph \fh_i = \ph\sum_{j=0}^m \vh_j [f_j,f_i]^x +{H}_{v_iu}\dot{\uh}=\hat{p}\, [f_0,f_i]^x,
\eenl
where Lemma \ref{NC} is used in the last equality. Therefore, 
$
\dot H_v = \dot H_v \lamh (\xh,\uh).
$
We can then differentiate $\dot{H}_v$ one more time, replace the occurrence of $\dot\uh$ by $\Gamma$ and obtain \eqref{ddotHv0} as it was desired.
See that \eqref{ddotHv0} together with the boundary conditions
\begin{gather}
 \label{HvT}  H_v\lamh (\wh_T)=0,\\
 \label{dotHv0} \dot H_v \lamh (\wh_0) = 0,
\end{gather}
guarantee the second identity in the stationarity condition \eqref{stationarity}.

\noindent {\it Notation:} Denote by (OS) the set of equations consisting of \eqref{stateeq}-\eqref{finaleq}, \eqref{costateeq}-\eqref{pT}, \eqref{Hu0},  \eqref{ddotHv0} and the boundary conditions \eqref{HvT}-\eqref{dotHv0}.

\begin{remark}
 Instead of \eqref{HvT}-\eqref{dotHv0}, we could choose another pair of endpoint conditions among the four possible ones: $H_{v,0}=0,$ $H_{v,T}=0,$ $\dot H_{v,0}=0$ and $\dot H_{v,T}=0,$ always including at least one of order zero. The choice we made will simplify the presentation of the result afterwards.
\end{remark}

Observe now that the derivative with respect to $(u,v)$ of the mapping
$
(w,\lambda) \mapsto
\begin{pmatrix}
 H_u\lam(w)
\\
-\ddot H_v\lam(w)
\end{pmatrix}
$
is  given by
\be 
\label{Jac}
\J:=
\begin{pmatrix}
 H_{uu} & H_{uv}
\\
-\ds\frac{\partial \ddot{H}_v}{\partial u} & 
-\ds\frac{\partial \ddot{H}_v}{\partial v}
\end{pmatrix}.
\ee
On the other hand, if \eqref{strongLC} and \eqref{strongenLC} are verified, $\J$ is definite positive along $(\wh,\hat\lambda)$ and, consequently, it is nonsingular. In this case we may write $\uh=U\lamh (\xh)$ and $\vh=V\lamh (\xh)$ from  \eqref{Hu0} and \eqref{ddotHv0}. Thus (OS) can be regarded as a TPBVP whenever the following hypothesis is verified.

\begin{assumption}
\label{assumptionLC}
$(\wh,\hat\lambda)$ satisfies \eqref{strongLC} and \eqref{strongenLC}.
\end{assumption}

Summing up we get the following result.
\begin{proposition}[Elimination of the control]
\label{ElimCont}
If $\wh$ is a smooth weak minimum verifying Assumptions \ref{lambdaunique} and \ref{assumptionLC}, then 
\benl
\uh=U\lamh (\xh), \quad \vh=V\lamh (\xh),
\eenl
with smooth functions $U$ and $V.$
\end{proposition}
\begin{remark}
When the linear and nonlinear controls are uncoupled, this elimination of the controls is much simpler. An example is shown in Oberle \cite{Obe90} where a nonlinear control variable can be eliminated by the stationarity of the pre-Hamiltonian, and the remaining problem has two uncoupled controls, one linear and one nonlinear. Another example is the one presented in Section \ref{SecEx}.
\end{remark}

\subsection{The algorithm}

The aim of this section is to present a numerical scheme to solve system (OS). In view of Proposition \ref{ElimCont} we can define the following mapping.
\begin{definition}
Let $\shoot:\cR^n\times\cR^{n+d_{\eta},*} =: {\rm D}(\mathcal{S})  \rightarrow \ \cR^{d_{\eta}}\times \cR^{2n+2m,*}$ be the \textit{shooting function} given by
\benl
\begin{array}{rl}
\begin{pmatrix} 
x_0,p_0,\beta 
\end{pmatrix}
 =:\nu
&\mapsto
\, \shoot (\nu):=
\begin{pmatrix}
\eta(x_0,x_T)\\
p_0 +D_{x_0}\ell\lam(x_0,x_T)\\
p_T -D_{x_T}\ell\lam(x_0,x_T)\\
H_v\lam(w_T)\\
\dot H_v (w_0)
\end{pmatrix},
\end{array}
\eenl
where $(x,p)$ is a solution of  \eqref{stateeq},\eqref{costateeq},\eqref{Hu0},\eqref{ddotHv0} 
with initial conditions $x_0$ and $p_0,$  $\lambda:=(p,\beta);$ and where the occurrences of $u$ and $v$ were replaced by $u=U[\lambda](x)$ and $v=V[\lambda](x).$
\end{definition}

Note that solving (OS) consists of finding $\hat\nu \in {\rm  D}(\mathcal{S})$ such that
\be\label{S=0}
\mathcal{S}(\hat\nu)=0.
\ee
Since the number of equations in \eqref{S=0} is greater than the number of unknowns, the Gauss-Newton method is a suitable approach to solve it. The \textit{shooting algorithm} we propose here consists of solving the equation \eqref{S=0} by the Gauss-Newton method.

\subsection{The Gauss-Newton Method}
This algorithm solves the equivalent least squares problem
\benl 
\min_{\nu \in {\rm  D}(\mathcal{S})} \left|\mathcal{S}\begin{pmatrix} \nu \end{pmatrix} \right|^2.
\eenl
At each iteration $k,$ given the approximate value $\nu^k,$ it looks for $\Delta^k$ that gives the minimum of the linear approximation of problem
\be\label{chap2minnorm}
\min_{\Delta\in {\rm  D}(\mathcal{S})} \left|\mathcal{S}(\nu^k)+
\mathcal{S}'(\nu^k)\Delta\right|^2.
\ee
Afterwards it updates
\benl
\nu^{k+1}\leftarrow \nu^k+\Delta^k.
\eenl
In order to solve the linear approximation of problem \eqref{chap2minnorm} at each iteration $k,$ we look for $\Delta^k$ in the kernel of the derivative of the objective function, i.e. $\Delta^k$ satisfying
\benl
\mathcal{S}'(\nu^k)\tras \mathcal{S}'(\nu^k)\Delta^k + \mathcal{S}'(\nu^k)\tras \mathcal{S}(\nu^k)=0.
\eenl
Hence, to compute direction $\Delta^k$ the matrix $\mathcal{S}'(\nu^k)\tras \mathcal{S}'(\nu^k)$ must be nonsingular. Thus, {Gauss-Newton method will be applicable provided that $\mathcal{S}'(\hat\nu)\tras \mathcal{S}'(\hat\nu)$ is invertible,} where $\hat\nu:=(\xh_0,\ph_0,\hat\beta).$ It follows easily that  $\mathcal{S}'(\hat\nu)\tras \mathcal{S}'(\hat\nu)$ is nonsingular if and only if $\mathcal{S}'(\hat\nu)$ is one-to-one.

 Furthermore, since the right hand-side of system \eqref{S=0} is zero, it can be proved that the Gauss-Newton algorithm converges locally quadratically if the function $\shoot$ has Lipschitz continuous derivative.
The latter holds true here given the regularity hypotheses on the data functions.
This convergence result is stated in the proposition below. See e.g. Fletcher \cite{Fle80} for a proof.
\begin{proposition}
\label{Conv}
  If $\mathcal{S}'(\hat\nu)$ is one-to-one then the shooting algorithm is locally quadratically convergent.
\end{proposition}

\section{Second order sufficient condition}
\label{SecSSC}

In this section we present a sufficient condition for optimality proved in \cite{Aro11}, and we state in  Section \ref{SecConv} afterwards that this condition guarantees the local quadratic convergence of the shooting algorithm proposed above.

Given $(\xb_0,\ub,\vb)\in \cR^n\times \U\times \V,$ consider 
the \textit{linearized state equation}
\begin{align} 
\label{lineareq} 
\dot{\xb}_t &= F_{x,t}\xb_t + F_{u,t}\ub_t + F_{v,t}\vb_t,\quad {\rm a.e.}\  {\rm on}\ [0,T],\\
\label{lineareq0}
\xb(0) &= \xb_0,
\end{align}
where $F_{x,t}$ refers to the partial derivative of $F$ with respect to $x,$ i.e. $D_x F_t;$ and equivalent notations hold for the other involved derivatives.
Take an element $\wb\in \W$ and define the second variation of the Lagrangian function
\benl 
\Omega (\wb) := \half D^2\L\lamh (\wh)\, \wb^2.
\eenl
It can be proved that $\Omega$ can be written as 
\benl
\begin{split}
\Omega &   (\xb,\ub,\vb)= \,
\half D^2\ell(\xh_0,\xh_T)(\xb_0,\xb_T)^2  + \intT \left[\half\xb\tras H_{xx} \xb \right. \\
& \left. + \ub\tras H_{ux}\xb +
\vb\tras H_{vx}\xb + \half\ub\tras H_{uu}\ub + \vb\tras H_{vu}\ub\right] \dtt.
\end{split}
\eenl
Note that this mapping $\Omega$ does not contain a quadratic term on $\vb$ since $H_{vv}\equiv 0.$ Hence, one cannot state a sufficient condition in terms of the uniform positivity of $\Omega$ on the set of critical directions, as it is done in the totally nonlinear case. Therefore, we use a change of variables introduced by Goh in \cite{Goh66a} and transform $\Omega$ into a quadratic mapping that may result uniformly positive in an associated transformed set of critical directions.

Consider hence the linear differential system in \eqref{lineareq} and
the change of variables 
\be \label{Goht}
\left\{
\ba{l}
\yb_t:= \ds\int_0^t \vb_s {\rm d}s, \\
\xib_t := \xb_t-  F_{v,t}\,\yb_t, 
\ea
\right.
\quad {\rm for}\ t\in [0,T].
\ee
This change of variables can be done in any linear system of differential equations, and it is often called \textit{Goh's transformation}.
Observe that $\xib$ defined in that way satisfies the linear equation
\be \label{xieq}
\dot\xib = F_x\xib + F_u\ub +B\yb,\quad
\xib_0=\xb_0,
\ee
where $ B:= F_{x}F_{v}-\ddt F_{v}.$

\subsection{Critical cones} 
We define now the sets of critical directions associated with $\wh.$ Even if we are working with control variables in $L^{\infty}$ and hence the control perturbations are naturally taken in $L^{\infty},$ the second order analysis involves quadratic mappings and it is useful to extend them continuously to $L^2.$ 
Given $\wb\in \W_2:= W^{1,2}(0,T;\cR^n) \times L^2(0,T;\cR^l) \times L^2(0,T;\cR^m)$ satisfying \eqref{lineareq}-\eqref{lineareq0}, consider
the \textit{linearization of the endpoint constraints and cost function,}
\begin{gather}
\label{linearconseq}
D\eta_j(\xh_0,\xh_T)(\xb_0,\xb_T)=0,\quad {\rm for}\  j=1,\hdots,d_{\eta},
\\
\label{linearconsineq}
 D\varphi_0(\xh_0,\xh_T)(\xb_0,\xb_T)\leq 0.
\end{gather}
Define the \textit{critical cone} in $\W_2$ by
\be
\C_2:=\{\wb\in\W_2:\text{\eqref{lineareq}-\eqref{lineareq0},\,\eqref{linearconseq}-\eqref{linearconsineq}}\ \text{hold}\}.
\ee
Since we aim to state an optimality condition in terms of the variables after Goh's transformation, we transform
the equations defining $\C_2.$   Let $(\xb,\ub,\vb)\in \C_2$ be a critical direction. Define $(\xib,\yb)$ by transformation
\eqref{Goht} and set $\hb:=\yb_T.$ Then the transformed of \eqref{linearconseq}-\eqref{linearconsineq} is 
\begin{gather}
\label{tlinearconseq}
D\eta_j (\xh_0,\xh_T)(\xib_0,\xib_T+B_T\hb)=0,\ 
\mr{for}\,\, j=1,\hdots,d_{\eta},
\\
\label{tlinearconsineq}
D\varphi_0(\xh_0,\xh_T)(\xib_0,\xib_T+B_T\hb)\leq
0.
\end{gather}
Consequently, the transformed critical cone is given by
\be
\P_2:= \{(\xib, \ub,\yb,\hb)\in \W_2\times
\cR^m:\,\text{\eqref{xieq}, \eqref{tlinearconseq}-\eqref{tlinearconsineq}
hold} \}.
\ee
\subsection{Second variation}
Next we state that if $\wh$ is a weak minimum, then the transformation of $\Omega$ yields the quadratic mapping
\be
\label{OmegaP2}
\ba{r}
\vspace{3pt}\bar\Omega(\xib,\ub,\yb,\hb):= g
(\xib_0,\xib_T,\hb)
+ \ds \intT \left( \half\xib\,\tras H_{xx}\xib + \ub\tras
H_{ux}\xib \right.\\
 \left. +\, \yb\tras M\lam \xib + \half\ub\tras H_{uu}\lam
\ub + \yb\tras J\lam \ub + \half\yb\tras R\lam \yb
\right) \dtt,
\ea
\ee
with 
\begin{gather*}
M:= F_v\tras H_{xx}-\dot H_{vx}-H_{vx}F_x,\ J:= F_v\tras H_{ux}\tras - H_{vx}F_u,\\
S:=\half (H_{vx}F_v+(H_{vx}F_v)\tras),\\
V:= \half (H_{vx}F_v-(H_{vx}F_v)\tras),
\\
R := F_v\tras H_{xx}F_v - (H_{vx}B+(H_{vx}B)\tras) -
\dot S,
\end{gather*}
\benl
g(\zeta_0,\zeta_T,h):=
\half\ell''(\zeta_0,\zeta_T+F_{v,T}h)^2
+h\tras(H_{vx,T} \zeta_T+\half S_T h).
\eenl
Easy computations show that $V_{ij}=\ph[f_j,f_i]^x,$ for $i,j=1,\dots,m.$ Thus, in view of Lemma \ref{NC}, one has that $V\equiv 0$ if $\wh$ is a weak minimum. Furthermore, we get the following result, which also uses \cite[Theorem 4.4]{Aro11}.

\begin{theorem}
If $\wh$ is a smooth weak minimum, then 
\benl
\Omega(\xb,\ub,\vb)=\bar\Omega(\xib,\ub,\yb,\yb_T),
\eenl
for all $(\xb,\ub,\vb)\in \W$ and $(\xib,\ub,\yb)$ given by  \eqref{Goht}.
\end{theorem}

\subsection{The sufficient condition}
We state now a second order sufficient
condition for strict weak optimality. 

Define the  \textit{$\gamma-$order} by
\benl 
\bar\gamma(\bar\zeta_0,\ub,\yb,\hb):=|\bar\zeta_0|^2+ |\hb|^2+
\intT(|\ub_t|^2+|\yb_t|^2)\dtt,
\eenl
for $(\bar\zeta_0,\ub,\yb,\hb)\in \cR^n\times L^2(0,T;\cR^l)\times  L^2(0,T;\cR^m) \times
\cR^{m}.$
It can also be considered as a function of
$(\bar\zeta_0,\ub,\vb)\in \cR^n\times L^2(0,T;\cR^l)\times L^2(0,T;\cR^m)$ by setting
\be 
\gamma(\bar\zeta_0,\ub,\vb):=
\bar\gamma(\bar\zeta_0,\ub,\yb,\yb_T),
\ee
with $\yb$ being the primitive of $\vb$ defined in
\eqref{Goht}.

\begin{definition}\label{qgdef}[$\gamma-$growth]
We say that $\wh$ satisfies
$\gamma-$\textit{growth condition in the
weak sense} if there exist $\eps,\rho > 0$ such
that 
\be 
\label{qg}
J(w) \geq J(\wh) + \rho \gamma(x_0-\xh_0,u-\uh,v-\vh),
\ee 
for every feasible trajectory $w$ with $\|w-\wh\|_{\infty} < \eps.$
 \end{definition}


\begin{theorem}[Sufficient condition for weak optimality]
\label{SC}
Let $\wh$ be a smooth feasible trajectory such that Assumption \ref{lambdaunique} is satisfied. Then the following assertions hold.
\begin{itemize}
\item[(i)]
Assume that there exists $\rho > 0$ such that
\be
\label{unifpos}
\bar\Omega (\xib,\ub,\yb,\hb) \geq
\rho\bar\gamma(\xib_0,\ub,\yb,\hb),\quad \text{on}\
\P_2.
\ee
Then $\wh$ is a weak minimum satisfying 
$\gamma-$growth in the weak sense.
\item[(ii)] 
Conversely, if $\wh$ is a weak solution satisfying $\gamma-$growth in the weak sense then \eqref{unifpos} holds for some $\rho>0.$
\end{itemize}
\end{theorem}

\section{Main result: Convergence of the shooting algorithm}\label{SecConv}
The main result of this article is the theorem below that gives a condition guaranteeing the quadratic convergence of the shooting method near an optimal local solution.
\begin{theorem}\label{wp} 
Suppose that $\wh$ is a smooth weak minimum satisfying Assumptions \ref{lambdaunique} and \ref{assumptionLC}, and such that \eqref{unifpos} holds. Then the shooting algorithm is locally quadratically convergent.
\end{theorem}

\begin{remark}
The complete proof of this theorem can be found in \cite{Aro11}. The idea of the proof is to show that \eqref{unifpos} yields the injectivity of $\shoot'(\hat\nu)$ and then use Proposition \ref{Conv}. In order to prove that \eqref{unifpos} implies that $\shoot'(\hat\nu)$ is one-to-one, the following elements are employed: the linearization of (OS) which gives an expression of the derivative $\shoot'(\hat\nu),$ the Goh's transformed of this linearized system and an associated linear-quadratic optimal control problem in the variables $(\xib,\ub,\vb,\hb)$ involving \eqref{xieq} and \eqref{OmegaP2}.
\end{remark}

\begin{remark}
[Bang-singular solutions] Finally we claim that the formulation of the shooting algorithm above and the proof of its local convergence (Theorem \ref{wp}) can be done also for problems where the controls are subject to bounds of the type
\be
\label{bound}
0\leq u_t \leq 1,\quad 0\leq v_t \leq 1,\quad \text{a.e.}\ \text{on}\ [0,1].
\ee
More precisely, it holds for solutions for which each control component is a  concatenation of {\it bang} and singular arcs, i.e. arcs saturating the corresponding inequality in \eqref{bound}, and arcs in the interior of the constraint.
This extension follows from a transformation of the problem to one without bounds, and it is detailed in  \cite[Section 8]{ABM11} for the totally-affine case.
\end{remark}

\if{

\section{The control constrained case: bang-singular solutions}\label{SecBounds}

In this section we add the following bounds to the control variables
\be
\label{chap2contcons}
0\leq u_{i,t}\leq 1,\quad {\rm  for}\ {\rm  a.a.}\ t\in [0,T],\ {\rm  for}\ i=1,\hdots,m.
\ee
Denote with (CP) the problem given by \eqref{chap2cost}-\eqref{chap2finalcons} and \eqref{chap2contcons}.
\begin{definition}
A feasible trajectory $\wh\in\W$ is  a \textit{Pontryagin minimum} of (CP) if for any positive $N$ there exists $\varepsilon_N>0$ such that $\wh$ is a minimum in the set of feasible trajectories $w=(x,u)\in\W$ satisfying 
\benl
\|x-\xh\|_{\infty}<\varepsilon_N,\ \|u-\uh\|_{1}<\varepsilon_N,\ \|u-\uh\|_{\infty}<N.
\eenl
\end{definition}

Given $i=1,\hdots,m,$ we say that $\uh_i$ has a \textit{bang arc} in $(a,b)\subset ]0,T[$ if $\uh_{i,t}=0$ a.e. on $(a,b)$ or $\uh_{i,t}=1$ a.e. on $(a,b),$ and it has a \textit{singular arc} if $0<  \uh_{i,t} < 1$ 
a.e. on $(a,b).$

\begin{assumption}
\label{chap2geohyp}
Each component $\uh_i$ is a finite concatenation of bang and singular arcs.
\end{assumption}

A time $t\in ]0,T[$ is called \textit{switching time} if there exists an index $1\leq i\leq m$ such that $\uh_i$ switches at time $t$  from singular to bang, or vice versa, or from one bound in \eqref{chap2contcons} to the other.

\begin{remark}
Assumption \ref{chap2geohyp} rules out the solutions having an infinite number of switchings in a bounded interval. This behavior is usually known as Fuller's phenomenon (see Fuller  \cite{Ful63}).
Many examples can be encountered satisfying Assumption \ref{chap2geohyp} as is the case of the three problems presented in Section \ref{chap2SecNum}.
\end{remark} 

With the purpose of solving (CP) numerically we assume that the structure of the concatenation of bang and singular arcs of the optimal solution $\wh$ and an approximation of its switching times are known. 
This initial guess can be obtained, for instance, by solving the nonlinear problem resulting from the discretization of the optimality conditions or by a continuation method. See Betts \cite{Bet98} or Biegler \cite{Bie10} for a detailed survey and description of numerical methods for nonlinear programming problems. For the continuation method the reader is referred to Martinon \cite{MartinonThesis}.

This section is organized as follows. From (CP) and the known structure of $\uh$ and its switching times we create a new problem that we denote by (TP). Afterwards we prove that we can transform $\wh$ into a weak solution $\Wh$ of (TP). Finally we conclude that if $\Wh$ satisfies the coercivity condition \eqref{chap2unifpos}, then the shooting method for problem (TP)  converges locally quadratically.
In practice, the procedure will be as follows: obtain somehow the structure of the optimal solution of (CP), create problem (TP), solve (TP) numerically obtaining $\Wh,$ and finally transform $\Wh$ to find $\wh.$

Next we present the transformed problem.

\begin{assumption}
\label{chap2disc}
 Assume that each time a control $\uh_i$ switches from bang to singular or vice versa, there is a discontinuity of first kind.
\end{assumption}

Here, by \textit{discontinuity of first kind} we mean that each component of $\uh$ has a finite nonzero jump at the switching times,
and the left and right limits exist.

By Assumption \ref{chap2geohyp} the set of switching times is finite. Consider the partition of $[0,T]$ induced by the switching times:
\be\label{chap2part}
\{0=:\Th_0< \Th_1<\hdots < \Th_{N-1}< \Th_N:=T\}.
\ee
Set $\hat{I}_k:=[\Th_{k-1},\Th_k],$ and define for $k=1,\hdots,N,$
\begin{align}
S_k&:= \{ 1\leq i\leq m:\ \uh_i\ {\rm  is}\ {\rm singular}\ {\rm  on}\ \hat{I}_k\}, \\
E_k&:= \{ 1\leq i\leq m:\ \uh_i=0\ {\rm  a.e.}\ {\rm  on}\ \hat{I}_k \},\\
N_k&:= 
\{ 1\leq i\leq m:\ \uh_i=1\ {\rm  a.e.}\ {\rm  on}\ \hat{I}_k\}.
\end{align}
Clearly $S_k \cup E_k \cup N_k=\{1,\hdots,m\}.$

\begin{assumption}
\label{chap2strongLC}
 For each $k=1,\hdots,N,$ denote by $u_{S_k}$ the vector with components $u_i$ with $i\in S_k.$ Assume that the strengthened generalized Legendre-Clebsch condition  holds on $\hat{I}_k,$ i.e. 
\be\label{chap2strongLCeq}
-\frac{\partial}{\partial u_{S_k}}\ddot{H}_{u_{S_k}} \succ 0,\quad {\rm  on}\ \Ih_k.
\ee
\end{assumption}

Hence,  $u_{S_k}$ can be retrieved from equation
\be
\label{chap2ddotHuSk}
\ddot{H}_{u_{S_k}}=0,
\ee
since the latter is affine on $u_{S_k}$
as it has been already pointed out in Section \ref{chap2SecOS}. 
Observe that the expression obtained from \eqref{chap2ddotHuSk} involves only the state variable $\xh$ and the corresponding adjoint state $\ph.$ Hence, it results that $\uh_{S_k}$ is continuous on $\Ih_k$ with finite limits at the endpoints of this interval.
 As the components $\uh_i$ with $i\notin S_k$ are either identically 1 or 0, we conclude that
\be\label{chap2contcont}
\uh\ \text{is continuous on}\ \Ih_k.
\ee

By Assumption \ref{chap2disc} and condition \eqref{chap2contcont} (derived from  Assumption \ref{chap2strongLC}) we get that there exists $\rho >0$ such that 
\be
\label{chap2contpos}
\rho<\uh_{i,t} <1-\rho,\quad \text{a.e. on}\ \Ih_k,\ 
\text{for}\ k=1,\hdots,N, \ i\in S_k.
\ee
Next we present a new control problem obtained in the following way. For each $k=1,\hdots,N,$ we perform the change of time variable that converts the interval $\Ih_k$ into $[0,1]$, afterwards we fix the bang control variables to their bounds and finally, we associate a free control variable to each index in $S_k.$
 More precisely, consider for $k=1,\hdots,N$  the control variables $u_i^k\in L_{\infty}(0,1;\cR),$ with $i\in S_k,$ and the state variables $x^k\in W^1_{\infty}(0,1;\cR^n).$ 
Let the constants $T_k\in \cR,$ for $k=1,\hdots,N-1,$ which will be considered as state variables of zero-dynamics.
Set $T_0:=0,$ $T_N:=T$ and define the problem on the interval $[0,1]$
\begin{align}
 & \label{chap2cost2} \varphi_0(x^1_0,x^N_1)\rightarrow\min,\\
 & \label{chap2stateeq2} \dot x^k=(T_k-T_{k-1})
\left(\sum_{i\in N_k\cup \{0\}}f_i(x^k)+\sum_{i\in S_k} u_i^{k}f_i(x^k)\right),\quad \ k=1,\hdots,N,\\
 & \label{chap2Teq} \dot T_k=0,\quad \ k=1,\hdots,N-1, \\
 & \label{chap2finalcons2} \eta(x^1_0,x^N_1)=0,\\
 & \label{chap2continuity} x^k_1=x^{k+1}_0,\quad k=1,\hdots,N-1.
\end{align}
Denote by (TP) the problem consisting of equations \eqref{chap2cost2}-\eqref{chap2continuity}.
The link between the original problem (CP) and the transformed one (TP) is given in Lemma \ref{chap2link} below.
Set for each $k=1,\hdots,N:$
\begin{align}
 \label{chap2xk}\xh^k_s &:= \xh(\Th_{k-1}+(\Th_k-\Th_{k-1})s),\quad {\rm for}\ s\in [0,1],\\
 \label{chap2uk}\uh_{i,s}^k &:= \uh_i(\Th_{k-1}+(\Th_k-\Th_{k-1})s),\quad {\rm  for}\ i\in S_k,\ {\rm a.a.}\ s\in [0,1].
\end{align}
Set 
\be\label{chap2Whdef}
\Wh:=((\xh^k)_{k=1}^N,(\uh_ i^k)_{k=1,i\in S_k}^N,(\Th_k)_{k=1}^{N-1}).
\ee
\begin{lemma}\label{chap2link}
If $\wh$ is a Pontryagin minimum of (CP), then $\Wh$  is a weak solution of (TP).
\end{lemma}
}\fi

\section{An example}\label{SecEx}

Consider the following optimal control problem treated in Dmitruk and Shishov \cite{DmiShi10}:
\be
\label{P}
\begin{split}
& J := -2x_{1,1} x_{2,1}+ x_{3,1} \to \min, \\
& \dot{x}_1 = x_2+u, \\
& \dot{x}_2 = v, \\
& \dot{x}_3 = x_1^2 + x_2^2 + 10x_2v+u^2, \\
& x_{1,0}=0, \quad x_{2,0}=0, \quad x_{3,0}=0.
\end{split}
\ee
Here, Assumption 1 holds since no final constraints are considered.
The pre-Hamiltonian function associated with \eqref{P} is, omitting arguments,
\benl
H=p_1(x_2+u)+p_2v+p_3(x_1^2 + x_2^2 + 10x_2v+u^2).
\eenl
We can easily deduce that $p_3 \equiv 1.$ 
The equations \eqref{Hu0} and \eqref{ddotHv0} for this problem give
\be
\label{exeqs}
H_u = p_1 + 2u,\quad \ddot{H}_v = -2v+2x_1,
\ee
and, therefore, Assumption 2 holds true. 
Agreeing with Proposition \ref{ElimCont}, the control can be eliminated from \eqref{exeqs}. This yields
\benl
u=-{p_1}/{2}, \quad v=x_1.
\eenl
We can then write the optimality system (OS) related to \eqref{P}. The state and costate equations are 
\be
\begin{split}
& \dot{x}_1 = x_2-{p_1}/{2}, \\
& \dot{x}_2 = x_1, \\
& \dot{x}_3 = x_1^2 + x_2^2 + 10x_2x_1+p_1^2/4, \\
& \dot{p}_1=-2x_1, \\
& \dot{p}_2 = -2x_2-10x_1-p_1,
\end{split}
\ee
where we do not include $p_3$ since it is constantly equal to 1. The boundary conditions are
\be
\label{exboundary}
\begin{split}
& x_{1,0}=0, \quad x_{2,0}=0, \quad x_{3,0}=0, \\
& p_{1,1}=-2x_{2,1},\quad p_{2,1}= -2x_{1,1}, \\
& H_{v,1} = p_{2,1}+10x_{2,1} =0,\\
& \dot{H}_{v,1}=-2x_{2,1}-p_{1,1}=0.
\end{split}
\ee
Observe that the last line in \eqref{exboundary} can be removed since it is implied by the first equation in the second line.
Here the shooting function is given by
\be
\label{Sexample}
\shoot:\cR^2 \to \cR^3,\,\, (p_{1,0},p_{2,0}) {\mapsto} \begin{pmatrix} p_{1,1} +2 x_{2,1} \\ p_{2,1} +2x_{1,1} \\ p_{2,1} + 10x_{2,1} \end{pmatrix}.
\ee

In \cite{DmiShi10} it was checked that the second order sufficient condition \eqref{unifpos} held for the control $(u \equiv 0, v\equiv 0).$ The solution associated with this control has $x_1=x_1=x_3=p_1=p_2=0.$ In view of Theorem \ref{SC}, we know that the shooting algorithm converges quadratically for appropriate initial values of $(p_{1,0},p_{2,0}).$

We solved \eqref{P} numerically  by applying the Gauss-Newton method to the equation $\shoot(p_{1,0},p_{2,0})=0,$ for  $\shoot$ defined in \eqref{Sexample}.  We used implicit Euler scheme for numerical integration of the differential equation. For arbitrary guesses of $(p_{1,0},p_{2,0}),$ the algorithm converged to $(0,0);$ in all the occasions. 
The tests were done with Scilab.

\section{CONCLUSIONS}\label{SecConclusion}

We investigated optimal control problems with systems that are affine in some components of the control variable and that have finitely many equality endpoint constraints. For a Mayer problem of this kind of system we proposed a numerical indirect method for approximating a weak solution. For qualified solutions, we proved that the local convergence of the method is guaranteed by a second order sufficient condition for optimality proved before by the author.

We presented an example, in which we showed how to eliminate the control by using the optimality conditions, proposed a shooting formulation and solved it numerically. The tests converged, as it was expected in view of the theoretical result.

\addtolength{\textheight}{-12cm}   





\section*{ACKNOWLEDGMENT}

Part of this work was done under the supervision of Fr\'ed\'eric Bonnans during my Ph.D. study. I acknowledge him for his great guidance.

I also thank Xavier Dupuis for his careful reading, and the three anonymous reviewers for their useful remarks.

This work is supported by the European Union under the 7th Framework Programme FP7-PEOPLE-2010-ITN  Grant agreement number 264735-SADCO.


\end{document}